\documentclass{amsart}

\usepackage{amssymb}
\newtheorem{thm}{Theorem}

\newtheorem{prop}{Proposition}

\begin{document}

\bibliographystyle{plain}

\title[Romi F. Shamoyan ]{On  an  extremal  problem  in  analytic  spaces  in  two  Siegel  domains  in  $C^{n}$}

\author[]{Romi F. Shamoyan}

\address{Department of Mathematics, Bryansk State Technical University, Bryansk ,241050, Russia}
\email{\rm rshamoyan@yahoo.com}

\date{}

\begin{abstract}
New sharp estimates concerning distance function in certain
Bergman -type spaces of analytic functions in a certain Siegel
domain of first type are obtained. Related sharp new estimates for
more general Siegel domains of  second type are also provided. For
Siegel domains of  second type in $C^{n}$ these are the first
results of this type.
\end{abstract}

\maketitle

\footnotetext[2]{\, Mathematics Subject Classification 2010
Primary 42B15, Secondary 42B30.  Key words and Phrases: Distance
estimates,analytic functions,Siegel domains of  first type and
second type}

\section{Introduction }

In this paper we obtain sharp distance estimates in certain spaces
of analytic functions in Siegel domains of  first type and of
second type. These types of domains are known in literature.They
studied by many authors during last decades (see for example
\cite{SG},\cite{BK},\cite{B1},\cite{SV} and references there). In
connection with the study of authomorphic functions of several
complex variables the notion of Siegel domains of the first and of
the second type was introduced by Piatetskii-Shapiro.(see\cite{SG}
and \cite{SV}).We recall basic facts which relate them to some
well-known domains.The Siegel domain of first type is a particular
case of a Siegel domain of second type(see \cite{SG})and in
particular there is a Siegel domain isomorphic to unit ball of
$C^{m+1}$ and in addition the simplest case of one dimensional
Siegel domain of the fist kind is a domain which we call usually
upperhalfspce $C_{+}$.Note also our results below were already
proved in this case  in \cite{SA}.Next the Siegel domain of first
type is a special type of a actively studied recently general tube
domains over symmetric cones (see \cite{SH} and various references
there concerning tube domains).But note there are homogeneous
Siegel domains of second type which are not even symmetric domains
(see \cite{SG} and \cite{SV}).Tube domains also are special cases
of Siegel domains of second type.It is known that every bounded
homogeneous domain in $C^{n}$ can be realized as Siegel domain of
the first and the second type and that this realization is unique
up to affine transformations. Siegel domains are holomorphically
equivalent to a bounded domains.But a lot of bounded domains are
not holomorphically isomorphic to Siegel domains.(see\cite{SG}) We
will provide definitions of  Siegel domains of first type and more
general of second type below referring also to \cite{SV} .(see
also, for example , \cite{SG})

Our line of investigation in this work can be considered as direct
continuation of our previous papers on extremal problems (see
\cite{AS1}, \cite{SM1} and \cite{SM2}). Our main two new results
are contained in the second and third sections of this note.First
we provide a concrete special example of a Siegel domain of first
type and we obtain a sharp estimate for distance function in
certain Bergman type analytic spaces on it.Next we turn in our
final section to Siegel domains of the second type. We remark that
here for the first time in literature we consider this extremal
problem related with distance estimates in spaces of analytic
functions in Siegel domains of second type. The next two sections
partially also contain some required preliminaries on analysis on
these domains.

In  the upperhalfspace  $C_{+}$  which is one dimensional tubular
domain and also in general tubular domains our theorems are not
new and they were obtained recently in \cite{SA},and then in
general form in \cite{SH}. Moreover arguments in proofs we
provided below are similar to those we have in previous cases and
hence our arguments sometimes will be sketchy below . The base of
proof is again the so-called Bergman reproducing formula, but in
Siegel domains.(see ,for example,\cite{BD},\cite{BK} for this
integral representation and it is applications). This paper first
deals with a concrete example of Siegel domain of first type and
based on some results from \cite{BD} we present a sharp result in
this direction. But then in the final part we turn to more general
situation (see \cite{BK} for notation which will be constantly
needed in this second part )and we obtain some related estimates
for distances there also. Note again some results from \cite{BK}
are crucial here in last section for us .

We now  shortly remind the history of this extremal problem.

After the appearance of \cite {Zh} various papers appeared where
arguments which can be seen in \cite{Zh} were extended and changed
in various directions \cite {SM1},\cite{SM2},\cite{AS1}.

In particular in mentioned papers various new results on distances for analytic function spaces in higher dimension (unit ball and polydisk) were obtained.
Namely new results  for large scales of analytic mixed norm spaces in higher dimension were proved.

Later several new sharp results for harmonic function spaces of
several variables in the unit ball and upperhalfplane of Euclidean
space were also obtained (see ,for example, \cite{AS1} and
references there) The classical Bergman representation formula in
various  domains serves as a base in all these papers in proofs of
main results .Recently (see\cite{SHA}) concrete analogues of our
theorems were proved also in some spaces of entire functions of
one and several variables.
 Various other extremal problems in analytic function spaces also
were considered before in various papers( see for
example\cite{AL},\cite{KH},\cite{RU},\cite{KS}).In those just
mentioned papers other results around this topic and some
applications of certain other extremal problems can be found also.
                    \section {New sharp estimates for distances in analytic Bergman -type spaces in
Siegel domains of first type}
 This section is devoted to one of the main results of this paper.
We remark our notes namely this one and \cite{SH} are first papers
with sharp results on extremal problems in higher dimension in
$C^{n}$ , namely in analytic function spaces in Siegel domains in
$C^{n}$. We now establish some notation from \cite{BD} which will
be needed for us. Let $\Omega\subset C^{n}$ be an open nonempty
set.Let $W(\Omega)$ be the set of all weights in  $\Omega$.We mean
by this a set of all Lebesgue measurable functions acting from
$\Omega$ to $R_{+}$.For each such $\gamma$ function let
$L^{2}(\Omega,\gamma)$ be the Hilbert space of all $f$ functions
from $\Omega$ to $C$ so that the quazinorm
$\int_{\Omega}|f(z)|^{2}\gamma(z)dm(z)$ is finite,where $dm(z)$ is
a Lebesgue measure on $\Omega$.By $A^{2}(\Omega)$ we denote  the
analytic subspace of this space but for so -called special (see
\cite{BD}) admissible $\gamma$ weights and with the same quazinorm
(see \cite{BD}).Note for these weights it is a closed subspace of
$L^{2}(\Omega)$. Next according to well-known Riesz representation
theorem there is a unique function that for all functions from
this  space  a certain integral representation holds with a
certain function called Bergman kernel which is from
$L^{2}(\Omega)$ (see \cite{BD}and references there).In certain
cases and our case is of them in higher dimension this function
called Bergman kernel can be explicitly written.This last fact
alone already opens a large way for various investigations in this
research area. In the present paper first we look at the family of
the following admissible  weights $\gamma_\alpha(\tau)$,
$$\gamma_\alpha(\tau)=(\Im\tau_{1}-|\tilde\tau|^{2})^{\alpha} $$
$\alpha > -1$ on the concrete  Siegel domain of the first
type.(see \cite{BD})
$$\Omega=\left\{\tau \in C^{n} ,   \Im \tau_{1} >|\tilde{\tau}|^{2}\right\}$$
here we denote by $\tau$ and $\tilde{\tau}$ the following vectors
$\tau=(\tau_{1},\ldots,\tau_{n})$,$\tilde{\tau}=(\tau_{2},\ldots,\tau_{n})$
Let $w$ be a vector from $C^{n}$. Let also $dm_\beta(w)=(\Im
w_{1}-|\tilde w|^{2})^{\beta}dm(w)$,where $dm(w)$ is a Lebesgue
measure on $R^{2n}$ and we also define a Bergman kernel as
see\cite {BD}
$$B_\beta(\tau,w)=({\tau-\overline w})^{-\beta-n-1}=(u-2v)^{-n-1-\beta}$$
$u=i(\overline{\tau_{1}}-w_{1})$,$v=({\tilde{w}}{\tilde{\tau}})$,where
the last expression is as usual a scalar product of two vectors in
$C^{n-1}$ .These definitions are crucial for our paper. The goal
of this section to develop further  some ideas from our recent
already mentioned papers and to present a new sharp theorem in
mentioned Siegel domain of first type .

For formulation of our result we will now need various standard
definitions from the theory of these Siegel domains of first
type.( see \cite{SG},\cite{BD},)

Let $\Omega$ be the Siegel domain . $\mathcal H(\Omega)$ denotes the space of all holomorphic functions on $\Omega$.
Let further,for all positive $\beta$.
\begin{equation}
A^\infty_\beta (\Omega) = \left\{ F \in {\mathcal H}(\Omega) : \|
F \|_{A^\infty_\beta} = \sup_{x+iy \in \Omega}|F(x+iy)|
\gamma_\beta (x+iy) < \infty \right\},
\end{equation}
(we use in this paper the following notation $w=u+iv$ and
$z=x+iy$,$w\in \Omega$,$z\in \Omega$). It can be checked that this
is a Banach space.

For $1 \leq p < + \infty$ , $\alpha>-1$ we denote by $A_\alpha^{p}(\Omega)$ the weighted

Bergman space consisting of analytic functions $f$ in $\Omega$ such that

$$ \| F \|_{A_\alpha^{p}} =  \left( \int_\Omega |F(z)|^p \gamma_\alpha (z)dm(z) \right)^{1/p} < \infty $$.This is a Banach space.
Below we will restrict ourselves to $p=2$ case following
\cite{BD}. Replacing above A by L we will get as usual the
corresponding larger space $L^{2}_\nu(\Omega)$ of all measurable
functions in our domain $\Omega$  with the same quazinorm (see
\cite{BD}). The (weighted) Bergman projection $P_\beta$ is the
orthogonal projection from the Hilbert space $L^2_\nu(\Omega)$
onto its closed subspace $A^2_\nu(\Omega)$ and it is given by the
following integral formula (see \cite{BD})

\begin{equation}\label{bpro}
P_\beta f(z) = C_\beta \int_{\Omega} B_\beta (z, w) f(w) d m_\beta
(w),
\end{equation}
where $C_\beta$ is a special  constant (see \cite{BD}) and $\beta>
\frac{\nu-1}{2}$.For these values of $\beta$ this is a linear
bounded operator from $L^{2}_\nu$ to $A^{2}_\nu$.Hence using these
facts  we have that for any analytic function from
$A^2_\nu(\Omega)$  the following integral formula  is valid for
all functions from $A^{2}_\nu$ ,for all
$\beta$,$\beta>\frac{\nu-1}{2}$ and $\nu>-1$ (see\cite{BD})

\begin{equation}\label{bpro1}
f(z)=C_\beta \int_{\Omega} B_\beta (z, w) f(w) d m_\beta (w),
\end{equation}

In this case  sometimes below we say simply that the analytic $f$
function allows Bergman representation via Bergman kernel with
$\beta$ index.

We need also the following estimate (A) of Bergman kernel from
\cite{BD}.Let $t>-1$ and $\beta >0$.Then there is a positive
constant $c=c_{n,t,\beta}$ so that
$$\int_{\Omega}\gamma_{t}(\tau)|B_{t+\beta}(\tau,w)|dm(\tau)\leq c \gamma^{-1}_{\beta}(w)$$,$w \in \Omega$
This estimate of Bergman kernel will be used and not once below during the proof of our first theorem .

Note here also these assertions we just mentioned have direct
analogues in simpler cases of analytic function spaces in unit
disk,polydisk,unit ball,upperhalfspace $C_{+}$ and in spaces of
harmonic functions in the unit ball or upperhalfspace of Euclidean
space $R^{n}$ .These classical facts are well- known and can be
found ,for example , in some items in references (see, for example
,\cite{Zh},\cite{DS}).

Above and throughout the paper  we write $C$( sometimes with indexes ) to denote   positive  constants which might be different  each time we see them
(and even in a chain of inequalities), but is independent of the functions or variables being discussed.

As in case of analytic functions in unit disk,polydisk,unit ball,
and upperhalfspace $C_{+}$,and tubular domains over symmetric
cones, and in case of spaces of harmonic functions in Euclidean
space \cite {Zh},\cite{SA},\cite{AS1},\cite{SM1}, \cite{SM2} the
role of the Bergman representation formula and  estimates for
Bergman kernel are crucial in these issues related with our
extremal problem and our proof will be heavily based on them.

And as it  was mentioned already above a variant of  Bergman
representation formula is available also in Bergman- type analytic
function spaces in Siegel domains and this known fact (see
\cite{BD},\cite{SV},\cite{BK}), which is crucial  in various
problems in analytic function spaces in Siegel domains of both
types  is also used in our proof below.

Moreover will also need for our proof the following  additional
facts on integral representation of functions on these $\Omega$
domains which follows from assertions we already formulated above.
Note first that for all functions from $A^\infty_\alpha$ the
integral representations of Bergman we mentioned above with
Bergman kernel
$$B_{\nu}(z,w)$$ (with $\nu$ index) is valid for large enough
$\nu$.This follows directly from the fact that $A^\infty_\alpha $
for any $\alpha$ is a subspace of $ A^{2}_\tau$ if $\tau$ is large
enough.(see \cite{BD}). Moreover it can be easily shown that we
have a continuous embedding $A^2_\alpha \hookrightarrow
A^\infty_\beta $( see ,for example, \cite{BD} where the proof can
be found also)for a concrete $\beta$ depending on $\alpha$,
$\alpha>-1$ and this naturally leads to a problem of estimating
$${\rm dist}_{A^\infty_{\beta}}(f, A^2_\alpha)$$ for a given $f\in
A^\infty_\beta$,where $\beta=\frac{\alpha+n+1}{2}$,$\alpha>-1$.

This problem on distances we just formulated will be solved in our
next theorem below, which is one of the main results of this
section. Let us set, for $f \in {\mathcal H}(\Omega)$, $s >0$ and
$\epsilon
> 0$ and  $z=x+iy \in \Omega $.
\begin{equation}\label{Efsep}
N_{\epsilon, s}(f) = \left\{ z \in \Omega : |f(z)|\gamma_{s}(z)
\geq \epsilon \right\}
\end{equation}
We denote by $N_{1}$ and by $N_{2}$ two sets- the first one is
$N_{\epsilon,s}(f)$, the other one is the set of all those points
,which are in tubular domain $\Omega$, but not in $N_{1}$ .Note
now ,to clarify the notation for readers again, by $m(z)$ or by
$m$ with only one lower index we denote in this section the
Lebesgue measure on $R^{2n}$

\begin{thm}\label{Td1}
Let $t = \frac{\nu+ 1+ {n}}{2}$.Set, for $f \in A^\infty_{\frac{n}{2} + \frac{\nu+1}{2}}$,$\nu >-1$

\begin{equation}\label{ETd1a}
l_1(f) = {\rm dist}_{A^\infty_{\frac{n}{2} + \frac{\nu+1}{2}}} (f, A^2_\nu),
\end{equation}
\begin{equation}\label{ETd1b}
l_2(f) = \inf \left\{ \epsilon > 0 : \int_{\Omega} \left(
\int_{N_{\epsilon, t}(f)} \frac{\gamma_{\beta - t }(w)
dw}{(z-\overline w )^{\beta+n+1}} \right)^2 \gamma_{\nu}(z) dm(z)
< \infty \right\}.
\end{equation}

Then there is a positive number $ \beta_{0}$   ,so that for all $ \beta > \beta_{0}$   we have $l_1(f) \asymp l_2(f)$.
\end{thm}

{\it Proof.}We will start the proof with the following
observation,which already was mentioned above. By our arguments
before formulation of this theorem for all functions from
$A^{\infty}_{\tau_{1}}$ the integral representations of Bergman
with Bergman kernel
$$B_(\tau_{2})(z,w)$$ is valid for large enough $\tau_{2}$

 We denote below the double integral which appeared in formulation of theorem by $G(f)$ and we will show first that  $l_{1}(f)\leq C l_{2}(f) $ .We assume now that $l_{2}(f)$ is finite.

 We use the Bergman representation formula which we provided above ,namely(\ref{bpro1}) ,and using conditions on parameters we now have  the following  equalities .

First we have obviously by remark from which we started this proof that for large enough $\beta$
$$f(z)=C_\beta\int_{\Omega}B_\beta(z,w)f(w)dm_\beta(w)= f_{1}(z)+f_{2}(z)$$
$$f_{1}(z)=C_\beta\int_{N_{2}}B_\beta(z,w)f(w)dm_\beta(w)$$,$$f_{2}(z)=C_\beta\int_{N_{1}}B_\beta(z,w)f(w)dm_\beta(w)$$

Then we estimate both functions separately using estimate (A)
provided above and following some arguments we provided in one
dimensional case that is the case of upperhalfspace $C_{+}$
\cite{SA}.Here our arguments are sketchy since they are parallel
to arguments from \cite{SA}. Using definitions of $N_{1}$ and
$N_{2}$ above  after some calculations  following arguments from
\cite {SA} using the estimate (A) of Bergman kernel we mentioned
above  we will have immediately.

$$f_{1}\in A^\infty_{\frac{\nu+n+1}{2}}$$   and  $$f_{2}\in A^2_\nu$$.We easily note the last inclusion follows directly from the fact that $l_{2}$ is finite.

Moreover it can be easily seen that the norm of $f_{1}$ can be
estimated from above by $C\epsilon$, for some positive constant
$C$ (\cite{SA}),since obviously
$$\sup_{N_{2}}|f(w)|\gamma_{t}(w)\leq \epsilon$$ .Note this last
fact follows directly from the definition of $N_{2}$ set and
estimate (A) above which leads to the following  inequality

$$\int_{\Omega}\gamma_{\beta-t}(\tau)|B_{\beta}(\tau,w)|dm(\tau)\leq C \gamma^{-1}_{t}(w)$$,$w \in \Omega$
for all $\beta$ so that $\beta>\beta_{0}$ ,for some large enough
fixed $\beta_{0}$ which depends on $n$,$\nu$, where
$$t=(\frac{1}{2})(\nu+1+n)$$

This gives immediately one part of our theorem .Indeed, we have
now obviously.

$$l_{1} \leq C_{2} \|f-f_{2}\|_{A^\infty_t} =C_{3} \|f_{1}\|_{A^\infty_t} \leq C_{4} {\epsilon}$$

It remains to prove that $l_{2} \leq l_{1}$.Let us assume $l_{1} <
l_{2}$.Then there are two numbers $\epsilon$ and
$\epsilon_{1}$,both positive such that there exists
$f_{\epsilon_{1}}$ ,so that this function is in $A^2_\nu$ and
$\epsilon > \epsilon_{1}$ and also the following  condition holds

$$\|f-f_{\epsilon_{1}}\|_{A^\infty_t}\leq \epsilon_{1}$$ and $G(f)=\infty$ ,where $G$ is a double integral  in formulation of theorem in $l_{2}$  .
(see (\ref{ETd1b}))

Next from $$\|f-f_{\epsilon_{1}}\|_{A^\infty_t} \leq
\epsilon_{1}$$ we have the following two estimates,the second one
is a direct corollary of first one.First we have

$$(\epsilon-\epsilon_{1})\tau_{N_{\epsilon,t}}(z)\gamma_{t}^{-1}(z)\leq C |f_{\epsilon_{1}}(z)|$$,where $\tau_{N_{\epsilon,t}}(z)$ is a characteristic function of $N=N_{\epsilon,t}(f)$ set we defined above.

And from last estimate we have directly multiplying both sides by
Bergman kernel $B_\beta(z,w)$ and integrating by tube $\Omega$
both sides with measure $dm_{\beta}$

$$G(f)\leq C \int_{\Omega}( L(f_{\epsilon_{1}}))^{2}\gamma_{\nu }(z) dm(z)$$, where $$L=L(f_{\epsilon_{1}},z)$$ and

$$L(f_{\epsilon_{1}},z)=\int_{\Omega}|f_{\epsilon_{1}}(w)||B_\beta(z,w)|dm_{\beta}(w)$$.Denote this expression by $I$.Put $\beta+n+1=k_{1}+k_{2}$,

where $k_{1}=\beta+1-n-\mu$,$k_{2}= \mu +
2n(\frac{1}{2}+\frac{1}{2})$ where the additional parameter will
be chosen by us  later.

By classical Holder inequality   we  obviously have

$$I^{2} \leq C I_{1} I_{2}$$, where

$$I_{1}(f)= \int_{\Omega}|f_{1}(z)|^{2}|(z-\overline w)^{s}|\gamma_{2\beta}(z)dm(z)$$

$$I_{2}= \int_{\Omega}|(z-\overline w)^{v}|dm(z)$$

and where $f_{1}=f_{\epsilon_{1}}$ and $$s=2\mu
-2-2\beta$$,$$v=-2n -2\mu $$.

Choosing finally the $\mu$ parameter ,so that the estimate (A)
namely

$$\int_{\Omega}\gamma_{\tilde{t}}(\tau)|B_{\tilde{t}+\tilde{\beta}}(\tau,w)|dm(\tau)\leq C \gamma^{-1}_{\tilde{\beta}}(w)$$,$$w \in
\Omega$$

 can be used twice above with some restrictions on parameters and finally making some additional easy calculations we will get what we need.

Indeed we have now obviously,

$$\int_{\Omega}(\int_{\Omega} |f_{\epsilon_{1}}(z)|B_{\beta}(z,w)|dm_{\beta}(z))^{2}\gamma_{\nu}(w)dm(w) \leq C
\|f_{\epsilon_{1}}\|^{2}_{A^{2}_{\nu}}$$ and

$$G(f)\leq C \|f_{\epsilon_{1}}\|_{A^2_\nu}$$,but we also have $$f_{\epsilon_{1}} \in A^2_\nu $$.

This will give as a contradiction with our previous assumption
above that $$G(f)=\infty$$. So we proved the estimate which we
wanted to prove. The proof of our first theorem in Siegel domains
of first type is now complete.
\section {New estimates for distances in Bergman type spaces in Siegel domains of the second
type}
          We first recall some basic facts on Siegel domains of second
type and then establish  notations for our second main theorem.
Recall first the explicit formula for the Bergman kernel function
is known for very few domains.The explicit forms and zeros of the
Bergman kernel function for Hartogs domains and Hartogs type
domains (Cartan-HArtogs domains) were found only
recently.\cite{AP1} .On the other hand in strictly pseudoconvex
domains the principle part of the Bergman kernel can be expressed
explicitly by kernels closely related to so-called Henkin -Ramirez
kernel see for example \cite{KH1} and references there. In
\cite{SG} the Bergman kernel
$$b\left((\tau_{1},\tau_{2}),(\tau_{3},\tau_{4})\right)$$ for the
Siegel domain of the second type was computed explicitly. It is an
integral via $V^{*}$  a convex homogeneous open irreducible cone
of rank $l$ in $R^{n}$ ,a conjugate cone of $V$ cone and which
also contains no straight line and in that integral the fixed
Hermitian form from definition of $D$ Siegel domain(see below for
definition) participates.(see details for this  \cite{BK}) .This
fact was heavily used in \cite{BK} in solutions of several
classical problems in Siegel domains of the second type and we
will also use one estimate from \cite{BK} for this kernel, but we
define it otherwise ,representing it otherwise in this paper.(see
also \cite{BK}).We will need now some short, but more concrete
review of certain results from \cite{BK} to make this exposition
more complete. To be more precise the authors in \cite{BK} showed
that on homogeneous Siegel domain of type 2 under certain
conditions on parameters the subspace of a weighted $L^{p}$ space
for all positive $p$ consisting of holomorphic functions is
reproduced by a concrete weighted Bergman kernel which we just
mentioned. They also obtain some $L^{p}$ estimates for weighted
Bergman projections in this case.The proof relies on direct
generalization of the Plancherel-Gindikin formula for the Bergman
space $A^{2}$ (see\cite{SG}).We remind the reader that the Siegel
domain of type 2 associated with the open convex homogeneous
irreducible  cone $V$ of rank $l$ which contains no straight
line,$V\in R^{n}$,and a $V$-Hermitian homogeneous form $F$ which
act from product of two $C^{m}$ into $C^{n}$ is a set of points
$(w,\tau)$ from $C^{m+n}$
 so that the difference $D$ of $\Im{w}$  and the value of $F$ on $(\tau,\tau) $ is in  $ V $ cone.This domain is affine homogeneous and we now should recall the following expression for the Bergman kernel
of  $$D=D(V,F)$$.Let $D$ be an affine-homogeneous Siegel domain of
type 2.Let $dv(z)$ denote the Lebesgue measure on $D$ and let
$$H(D)$$ denote the space of all holomorphic functions on $D$.The
Bergman kernel is given by the following formula( see\cite{BK})for
$(\tau_{1},\tau_{2})\in D$ and $(\tau_{3},\tau_{4})\in D$
$$b\left((\tau_{1},\tau_{2}),(\tau_{3},\tau_{4})\right)=(\frac{\tau_{1}-\overline\tau_{3}}{2i}-(F(\tau_{2},\tau_{4}))^{2d-q}$$ ,where two vectors $q=(q_{i})$ and
$d=(d_{i})$and in addition $n=(n_{i})$ here the $i$ index is
running from $1$ to $l$ are specified  via $n_{i,k}$ ,where these
 $n_{i,k}$ numbers are dimensions of certain $(R_{i,k})$ and
$(C_{i,j})$ subspaces of the certain canonical decomposition of
$C^{m+n}$ and $R^{n}$ via the $V$ cone from definition of our $D$
domain(see for some additional details about this \cite{SG} and
,\cite{BK}).We will call this family of triples parameters of a
Siegel domain $D$ of second type.They will appear in our main
theorem and it is short proof. The standard Bergman projection
here on $D$ as usual is denoted by $P$, it is the orthogonal
projection of $L^{2}(D,dv)$ onto it is analytic subspace
$A^{2}(D,dv)$ consisting of all holomorphic functions.The authors
in \cite{BK} showed that some well-known facts of much simpler
domains holds also here, for example there is an integral operator
on $L^{2}$space defined by the certain
$$b(\tau,z)$$ Bergman kernel.And for this types of  Siegel domains as it was
mentioned above this Bergman kernel was computed explicitly
previously in \cite{SG}.Further ,let $\epsilon$ be a real number.
Now for all positive finite $p$ we define a space of integrable
functions (weighted $L^{p}$ spaces with $b^{-\epsilon}(z,z)$
weights) for all $\epsilon>\epsilon_{0}$
$$L^{p,\epsilon}(D)=L^{p}(D,b^{-\epsilon}(z,z)dv(z))$$ and we
denote as usual by $A^{p,\epsilon}$ the analytic subspace of this
space with usual modification when $p=\infty$. Note the
restriction is meaningful since  there is an $\epsilon_{0}$ so
that for all those $\epsilon$ which are smaller than this fixed
$\epsilon_{0}$ the $A^{2,\epsilon}$ is an empty class (see
\cite{BK}). We denote by $P_{\epsilon}$ the corresponding Bergman
projection which is the orthogonal projection  of $L^{2,\epsilon}$
to it is analytic subspace $A^{2,\epsilon}$. In \cite{BK} the
authors give a condition on real numbers and vectors
$r,p,\epsilon$,  so that the weighted Bergman projection
reproduces all functions in $A^{p,r}(D)$.This vital fact for our
theorem they deduce partially from Plancherel-Gindikin formula and
the fact that
$$P_{\epsilon}(f)(z)=c_{\epsilon}\int_{D}f(u)b^{1+\epsilon}(z,u)b^{-\epsilon}(u)du$$ so it defines as in simpler cases an integral operator on
$L^{2,\epsilon}(D)$ by the kernel
$$b^{1+\epsilon}(\tau,z)$$ (see for this \cite{BK}), it is a
weighted Bergman projection from $L^{2,\epsilon}$ onto
$A^{2,\epsilon}$ (see, for example, \cite{BK} and references
there.) .The following several assertions concerning Bergman
projection acting in analytic spaces in Siegel domain of the
second type and estimates of Bergman kernel  which we mentioned
above and in addition to this some facts on spaces of integrable
functions and their analytic subspaces we defined above on these
Siegel domains were proved in \cite{BK} and some are crucial for
this paper. We will formulate immediately after them  our main
result on distances in Siegel domains of the second type. Then
providing a comment on a proof of that assertion which contains no
new ideas when we compare it with the proof of previous theorem we
will finish this paper.We use the following notation. The $i$
index below is running from $1$ to $l$ everywhere and to make the
reading easier we accept this from advance. We also use below
everywhere standard rules of calculations  between two vectors as
they were seen by us for example in \cite{BK},also sometimes we
write
$$d\tilde{V}(\tau_{1},\tau_{2})$$ not $dv(\tau)$ meaning
$$\tau=(\tau_{1},\tau_{2})\in D$$. In the following assertions
$$(n_{i}),(q_{i}),(d_{i})$$ are always act as parameters of the Siegel
domain $D$ we introduced above and they are  playing a crucial
role.We write always $D$ below  meaning $$D(n,q,d)$$,where
$n=(n_{i})$, $d=(d_{i})$,$q=(q_{i})$.We write $c_{i}\leq  b_{i}$
for two vectors from $R^{l}$ below meaning as usual that this is
true for all values of $i$ from $1$ to $l$.If $c\leq b{i}$ (or
$c<b_{i}$)then all $b_{i}$ are bigger or equal (or bigger)than
$c$.

\begin{prop}\label{THD}

Let $\epsilon\in R^{l}$, $r\in R^{l}$,$p\in R_{+}$,$0\leq
r_{j}$.Then there are two sets of numbers $(k_{i})$,$(m_{i})$
depending on parameters of $D$ Siegel domain so that if $1\leq p <
k_{i}$

and  $\epsilon_{i}> m_{i}$,then $$P_{\epsilon}f=f$$ for all
$$f\in A^{p,r}(D)$$

Let $\epsilon \in R^{l}$, $r\in R^{l}$,$p\in (0,\infty)$,$v_{i}<
r_{i}$,for some $v_{i}$ numbers depending on parameters of $D$
domain.Then there are two sets of numbers
$(k^{1}_{i})$,$(m^{1}_{i})$, depending from parameters of $D$
Siegel domain ,so that if $$0< p < k^{1}_{i}$$

and if $\epsilon_{i}>m^{1}_{i}$,then $$P_{\epsilon}f=f$$ for all
$$f\in A^{p,r}(D)$$

\end{prop}

\begin{prop}\label{THD1}

If $$\epsilon_{i}>\frac{n+2}{2(2d-q)_{i}}$$ where ,$\epsilon\in
R^{l}$ ,then $P_{\epsilon}$ is an integral operator with
$$b^{1+\epsilon}(t_{1},z_{1})(t_{2},z_{2})$$ kernel on
$L^{2,\epsilon}$ and $$P_{\epsilon}f=f$$  for all $f\in
A^{p,0}(D)$,when $p\in(0,p_{0})$,where $$p_{0}\leq
\frac{n_{i}-2(2d-q)_{i}}{n_{i}}$$.

If there is an index $i$ so that $$2\epsilon_{i}\leq
\frac{n_{i}+2}{(2d-q)_{i}}$$ then we have $A^{2,\epsilon}={0}$,
moreover if the reverse estimate holds for all $i$ and
$\epsilon_{i}$ instead of $2\epsilon_{i}$ then the intersection of
$A^{2,\epsilon}$ and $A^{p,r}$ is dense in $A^{p,r}$ ,if $1\leq
p<\infty$,$0\leq r_{i}$,$\epsilon\in R^{l}$

\end{prop}
The following embedding which is taken also from \cite{BK} is
important for us.It allows us as in previous simpler case to pose
a distance problem in this domain showing that  Bergman spaces
$A^{p,r}$ are subspaces of $A^{\infty}_{\frac{1+r}{p}}$
Bergman-type spaces. Let $r\in R^{l}$ and $p\in(0,\infty)$, then
$$|f(z)|^{p}\leq C b^{1+r}(z,z)\|f\|^{p}_{p,r}$$ $$z\in D$$.Further let
$\epsilon$ and $r$ are from $R^{l}$ .If
$$\epsilon_{i}>\frac{n_{i}}{-2(2d-q)_{i}}$$ and
$$r_{i}>\frac{n_{i}+2}{2(2d-q)_{i}}+\epsilon_{i}$$ then we have
$$P_{r}(f)=f$$ as soon as $f$ belongs to
$A^{\infty}_{\epsilon}$(see\cite{BK},\cite{AK}). This will also be
needed in proof of main result of this section.(see for this also
the parallel proof of our previous theorem from previous section)
\begin{prop}\label{THD2}

Let $\beta\in R^{l}$ and all $\beta_{i}$ are nonnegative then the
following estimate holds
$$b^{\beta}(z+\tau,z+\tau)\leq b^{\beta}(z,z)$$ and also $$|b^{\beta}(\tau,z)|\leq b^{\beta}(z,z)$$ for all $\tau$ and $z$ from $D$.
\end{prop}
The following estimate to be more precise it is direct analogue
can be found in the proof of previous theorem where it was used
three times.
\begin{prop}\label{THD3}
Let $\alpha$ and $\epsilon$ be two vectors from $R^{l}$ and
$(\tau,z)$ be a point from $D$ ,then if  $$
\frac{n_{i}+2}{2(2d-q)_{i}} < \epsilon_{i}$$ and if $$
\epsilon_{i}-\frac{n_{i}}{2(2d-q)_{i}} < \alpha_{i}$$ then the
following integral
$$\int_{D}|b^{\alpha+1}((\tau,v),(z,u))|b^{-\epsilon}((z,u)(z,u))d\tilde{V}(z,v)$$ is  equal to $$c_{\alpha,\epsilon}b^{\alpha-\epsilon}((\tau,v),(\tau,v))$$

\end{prop}
We are able now based only on last proposition and two comments
concerning integral representations before previous proposition to
formulate a theorem on distances in Siegel domains of the second
type which is a direct analogue of our previous results (see ,for
example, \cite{SM1},\cite{SM2},\cite{SA}) and our previous theorem
on distances in this  situation.All facts and preliminaries which
are needed here for our proof can be found above  in assertions
from \cite{BK} which we just formulated ,all lines of arguments
for our proof of this theorem can be also found above in the proof
of our previous theorem though some not very long technical
calculations with indexes should be added.Note one implication in
this theorem below is easier and we just repeat here arguments of
our previous theorem.
 \begin{thm}\label{Thd4}    Let  $$N_{\tilde\epsilon,r}(f)=\left\{z\in D ,|f(z)|b^{1+r}(z,z)>\tilde\epsilon \right\}$$,where $\tilde\epsilon $ is a positive
 number.Then the following two quantities are
 equivalent
 $${\rm dist}_{A^{\infty}_{1+r}}(f,A^{1,r})$$ and
 $$\inf\left\{\tilde\epsilon > 0,\int_{D}(\int_{N_{r,\tilde\epsilon}(f)}b^{-k+1+r}(\tau,\tau)|b(\tau,z)|^{k+1}dv(\tau))b^{-r}(z,z)dv(z)<\infty
 \right\}$$,for all $r$ and $k$ so that $r\in(r_{0},\infty)$ and $k\in(k_{0},\infty)$ and for certain
 fixed vectors $r_{0}$ and $k_{0}$   depending on parameters of the Siegel $D$ domain $(d_{i})$ and $(q_{i})$
 and $(n_{i})$
\end{thm}
We remark finally the  theorem above is probably valid for all
$p>1$ (not only $p=1$ when calculations are  simpler ) and the
reader can formulate easily that theorem in general case following
the formulation of our previous theorem .The proof probably is
parallel to the proof of previous theorem and it is based on
estimates from propositions above .Note also our assertion is true
for all homogeneous Siegel domains not only symmetric Siegel
domains of the second type(see \cite{BK},\cite{DB1},\cite{AK1}).
We remark as $r_{0}$ we can take $\max{({r_{1}},{r_{2}},{0})}$
where $r_{1}$ and $r_{2}$ are depending on parameters of domain
$r_{1}=\frac{n_{i}+2}{2(2d-q)_{i}}$ and
$r_{2}=\frac{-n_{i}}{2(2d-q)_{i}}-1$

\end{document}